\documentclass[sigplan,screen,nonacm]{acmart}
\usepackage{amsmath,bm}
\AtBeginDocument{%
  \providecommand\BibTeX{{%
    \normalfont B\kern-0.5em{\scshape i\kern-0.25em b}\kern-0.8em\TeX}}}

\setcopyright{acmcopyright}
\copyrightyear{2021}
\acmYear{2021}
\acmDOI{10.1145/1122445.1122456}

\acmConference{Preprint for ACM SIGEnergy Energy Informatics Review}  
\acmPrice{15.00}
\acmISBN{978-1-4503-XXXX-X/18/06}

\usepackage{color}

\usepackage{mathrsfs}
\usepackage{graphicx} 
\usepackage{subcaption}
\usepackage{float}
\usepackage[english]{babel}

\usepackage{tcolorbox}
\usepackage{enumerate}
\usepackage{bm}
\usepackage{url}
\usepackage{algorithm} 
\usepackage{algorithmic} 
\usepackage{mathtools}
\usepackage[mathcal]{euscript}

\DeclareMathOperator*{\argmin}{arg\,min}
\DeclareMathOperator{\Reg}{\texttt{Reg}}
\DeclareMathOperator{\Sreg}{\texttt{SReg}}
\DeclareMathOperator{\Cacv}{\texttt{CACV}}

\DeclareMathOperator{\st}{s. t.}

\newcommand{\BIT}{\begin{itemize}}
\newcommand{\EIT}{\end{itemize}}

\newcommand{\reals}{\mathbb{R}}

\newcommand{\bfx}{\mathbf{x}}

\newcommand{\calX}{\mathcal{X}}

\newcommand{\calV}{\mathcal{V}}
\newcommand{\calE}{\mathcal{E}}

\newcommand{\bfA}{\mathbf{A}}

\newcommand{\norm}[1]{\left\lVert#1\right\rVert}

\usepackage{color}
\newcommand*{\mathcolor}{}
\def\mathcolor#1#{\mathcoloraux{#1}}
\newcommand*{\mathcoloraux}[3]{%
  \protect\leavevmode
  \begingroup
    \color#1{#2}#3%
  \endgroup
}
\newcommand{\sts}{\mathcolor{white}{\st}}

\newenvironment{mybox}[1]{%
\tcolorbox[savedelimiter=mybox,
savelowerto=\jobname_bspsave2.tex,
lowerbox=ignored,
colback=red!5,colframe=red!75!black,fonttitle=\bfseries,title=#1]}%
{\endtcolorbox}

\begin{document}

\title{Towards   Online  Optimization for Power Grids}


\author{Deming Yuan}
\email{dmyuan1012@gmail.com}
\orcid{}
\affiliation{%
  \institution{School of Automation\\ Nanjing University of Science and Technology, Jiangsu 210094, China}
  \state{}
  \country{}
}

\author{Abhishek Bhardwaj}

\email{Abhishek.Bhardwaj@anu.edu.au}
\orcid{}
\affiliation{%
  \institution{School of Engineering \\ Australian National University}
  \state{ACT 2600}
  \country{Australia}
}

\author{Ian Petersen}
\email{ian.petersen@anu.edu.au}
\orcid{}
\affiliation{%
  \institution{School of Engineering \\ Australian National University}
  \state{ACT 2600}
  \country{Australia}
}


\author{Elizabeth L. Ratnam}
\email{elizabeth.ratnam@anu.edu.au}
\orcid{}
\authornotemark[1]
\affiliation{%
  \institution{School of Engineering \\ Australian National University}
  \state{ACT 2600}
  \country{Australia}
}

\author{Guodong Shi}
\email{guodong.shi@sydney.edu.au}
\orcid{}
\affiliation{%
  \institution{Australian Centre for Field Robotics\\ The University of Sydney}
  \state{NSW 2006}
  \country{Australia}
}
\authornote{Correspondence authors: E. L. Ratnam +61 429 369 924; G. Shi +61 2 8627 8037.}

\renewcommand{\shortauthors}{Yuan and Bhardwaj et al.}

\begin{abstract}
In this note, we discuss potential advantages in extending distributed optimization frameworks to enhance support for power grid operators managing an influx of online sequential decisions. First, we review the state-of-the-art distributed optimization frameworks for electric power systems, and explain how distributed algorithms  deliver scalable solutions. Next, we introduce key concepts and  paradigms for online optimization, and present a distributed online optimization framework highlighting important performance characteristics. Finally, we discuss the connection and difference between offline and online distributed optimization,  showcasing the suitability of such optimization techniques for power grid applications. 
\end{abstract}




\maketitle

 \section{Distributed Optimization for Power Grids}

\subsection{The evolving power grid}
The rise of distributed and renewable energy resources including wind and solar backed by energy storage technologies, is accelerating the evolution of the electric power grid \cite{Quint2019Transformation}. The evolution is supported by recent advances in communication \cite{Gungor2011Smart}, sensor technologies \cite{Meier2017Precision}, data processing \cite{Andersen2016BTrDB}, and operational technologies \cite{Erseghe2014Distributed, Yang2013Consensus, Elsayed2015Fully, Anese2013Distributed}, enabling prosumers to generate and deliver surplus renewable energy back to the power grid \cite{Chapman2021Network}. As the complexity of operating the evolving power grid continues to increase, with renewable technologies becoming increasingly distributed and spatially diverse, scalable approaches are needed to manage electricity flows to and from millions of energy prosumers. Distributed optimization frameworks that support scalability in managing renewable energy flows in transmission \cite{Huebner2019Distributed} and distributions networks \cite{Peng2018Distributed,Joo2017Distributed}, that additionally supporting solutions to integrate energy storage \cite{mcilwaine2021state, pimm2018potential}, including electrical vehicles \cite{NIMALSIRI2021Coordinated2, ghavami2016decentralized, NIMALSIRI2021Coordinated1}, potentially enhance both the operation and resilience of the evolving power grid.

\subsection{Distributed optimization frameworks}
The key promise in distributed optimization is a dramatic improvement in scalability to accommodate data and decisions scattered in physically decentralized locations. See \cite{molzahn2017survey} for an in-depth survey. 

\vspace{2mm} 

\noindent{\em Example 1.} (Economic Power Dispatch \cite{Elsayed2015Fully}) Consider $N$ generators indexed in $\mathcal{V}=\{1,\dots,N\}$. At a fixed time there is a total power demand $P$ that needs to be met by these $N$ generators. Let generator $i$ be allocated  a power $\mathbf{x}_i \in[x_i^{\rm min}, x_i^{\rm max}]$, leading to a cost $\ell_i(\mathbf{x}_i)$. The economic power dispatch problem:
\begin{align}
    \begin{split}
        \min_{\mathbf{x}} & \  \sum_{i=1}^{N} \ell_i(\mathbf{x}_{i}) \\
        \st & \ \ x_i^{\rm min}\leq \mathbf{x}_{i} \leq x_i^{\rm max}, \  i=1,\dotsc,N \\
        \sts & \ \ \mathbf{x}_{1}+\dotsb+\mathbf{x}_{N} = P. 
    \end{split} \label{ex1}
\end{align}
Here $\ell_i(\cdot)$ is a function mapping from $\mathbb{R}^{\geq 0}$ to $\mathbb{R}^{\geq 0}$. An optimal decision on the $\mathbf{x}_i$ for all $i$ should minimize the total generation cost.

\vspace{2mm} 

\noindent{\em Example 2.} (Optimal Power Flow \cite{Erseghe2014Distributed}) Consider an electrical network with $N$ nodes indexed in $\mathcal{V}=\{1,\dots,N\}$. Let $\mathbf{v}_i\in\mathbb{C}$ and $\mathbf{i}_i \in\mathbb{C}$ be the voltage and inflow current at node $i$. The network structure is captured by an admittance matrix $\mathbf{A}\in \mathbb{C}^{n\times n}$. Then $\mathbf{x}_i:={\rm Re} (\mathbf{v}_i \mathbf{i}_i^\dag)$ defines the active power at node $i$, where $^\dag$ is the complex conjugate. Let $\ell_i(\mathbf{x}_i)$ denote the cost associated with the power at node $i$. An optimal power flow problem is given in the following form:
\begin{align}
    \begin{split}
        \min_{\mathbf{x}} & \  \sum_{i=1}^{N} \ell_i(\mathbf{x}_{i}) \\
        \st &\ \  \mathbf{x}_{i}={\rm Re} \left(\mathbf{v}_i \mathbf{i}_i^\dag\right), \  i=1,\dotsc,N \\
        \sts &\ \   \mathbf{v}_i \mathbf{i}_i^\dag= \mathbf{v}_i\sum_{j=1}^n \mathbf{A}_{ij}^\dag \mathbf{v}_j^\dag, \  i=1,\dotsc,N.
    \end{split}\label{ex2}
    \end{align}

In centralized optimization for (\ref{ex1}) and (\ref{ex2}), each local cost function $\ell_i(\cdot)$ and local parameter such as $x_i^{\rm min},x_i^{\rm max},\mathbf{A}_{ij}$ needs to be sent a central coordinator; the coordinator solves the respective problem (\ref{ex1}) or (\ref{ex2}) and sends the optimal decisions for $\mathbf{x}_i$ to each agent. In distributed optimization for (\ref{ex1}) and (\ref{ex2}), there is an underlying communication graph $\mathcal{G}=(\mathcal{V},\mathcal{E})$ over which agents share their decisions, and computations are carried out locally in parallel at each individual agent based on the local cost functions and parameters. The distributed computing architecture naturally allows scalability; the absence of a central coordinator  improves resilience since failures at the coordinator have system-level impact while failures at individual agents harm the system-level performance at a limited level.

\subsection{Distributed optimization algorithms}
There are many algorithms for distributed optimization. In power systems, the Alternating Direction Method of Multipliers (ADMM) has been popular. 

\subsubsection{The ADMM}
Given an optimization problem 
\begin{align}\label{eqadmm}
   \begin{split}
 \min_{\mathbf{x},\mathbf{z}} & \ f(\mathbf{x}) +g(\mathbf{z})  \\
 \st & \ \mathbf{A}\mathbf{x}+\mathbf{B}\mathbf{z}=\mathbf{c},\\
 &\ \mathbf{x}\in\mathbb{R}^n, \mathbf{z}\in \mathbb{R}^m.  
    \end{split}
\end{align}
ADMM proceeds by first defining the augmented Lagrangian 
$$
L_{\alpha}(\mathbf{x},\mathbf{z},\mathbf{y}) = f(\mathbf{x})+g(\mathbf{z})+\mathbf{y}^{T}(\mathbf{A}\mathbf{x}+\mathbf{B}\mathbf{z}-\mathbf{c}) + \frac{\alpha}{2}\norm{\mathbf{A}\mathbf{x}+\mathbf{B}\mathbf{z}-\mathbf{c}}^{2}
$$
with dual variable $\mathbf{y}$. Then the algorithm runs recursively, where in each round there are updates in the decision variable for $\mathbf{x},\mathbf{z},\mathbf{y}$ that are arranged sequentially. 

\begin{mybox}{ADMM Algorithm \cite{2200000016}}
Define an initial point $(\mathbf{x}^{(0)}, \mathbf{z}^{(0)}, \mathbf{y}^{(0)})$, smoothing parameter $\alpha$, and iteration limit $n$.

For $k=0,\dotsc,n$ DO
\begin{itemize}
\item[(i)] Update the first variable as $$\mathbf{x}^{(k+1)} = \argmin_{\mathbf{x} \in\mathbb{R}^n} \  L_{\alpha}\left(\mathbf{x},\mathbf{z}^{(k)},\mathbf{y}^{(k)}\right) ;$$
\item[(ii)] Update the second variable as $$\mathbf{z}^{(k+1)} = \argmin_{\mathbf{z} \in\mathbb{R}^m} \  L_{\alpha}\left(\mathbf{x}^{(k+1)},\mathbf{z},\mathbf{y}^{(k)}\right) ;$$
\item[(iii)] Update the dual variable as $$\mathbf{y}^{(k+1)}= \mathbf{y}^{(k)} - \alpha \left( \mathbf{A}\mathbf{x}^{(k+1)}+\mathbf{B}\mathbf{z}^{(k+1)}-\mathbf{c}\right).$$
\end{itemize}
\tcblower

\end{mybox}

\subsubsection{Distributed ADMM}
The original ADMM algorithm was proposed in the 1970s \cite{Gabay1976ADA,M2AN}, and regained its popularity in recent years due to its suitability for large-scale distributed computing problems  \cite{2200000016}. 
%
If we write $f(\mathbf{x})=\sum_{i=1}^N \ell_i(\mathbf{x}_i)$ for the cost functions  in (\ref{ex1}) and (\ref{ex2}) with $\mathbf{x}=(\mathbf{x}_1,\dots,\mathbf{x}_N)$, and suitable auxiliary decision variables $\mathbf{z}$ from the constraints, problems in the form of (\ref{ex1}) and (\ref{ex2}) can be written in the standard ADMM form (\ref{eqadmm}). For example, the problem (\ref{ex1}) with $x_i^{\rm min}=-\infty$ and $x_i^{\rm max}=\infty$ can be written as (Chapter 7, \cite{2200000016})
\begin{align}\label{eqadmm1}
\begin{split}
 \min_{\mathbf{x},\mathbf{z}} & \ f(\mathbf{x}) +g_P\left(\sum_{i=1}^N \mathbf{z}_i\right)  \\
 \st & \ \mathbf{x}_i-\mathbf{z}_i=0, \ i=1,\dots,N,
 \end{split}
\end{align}
where $g_P(a)=1$ if $a=P$ and $g_P(a)=+\infty$ otherwise.  Then, due to the separable nature of the function $f$ and the constraints in (\ref{eqadmm1}), the resulting ADMM algorithm can be naturally decomposed into parallel computations at the agents along each primal variable $\mathbf{x}_i$ and dual variable $\mathbf{y}_i$ for Step (i) and Step (iii). The Step (ii) of the ADMM algorithm relies on all $\mathbf{x}_i$ and  $\mathbf{y}_i$ for $i=1,\dots,N$, and can be   implemented in a distributed fashion with the help of the communication graph $\mathcal{G}$.

\subsubsection{Alternatives to ADMM}

There are  many other distributed optimization methods aside from ADMM. Some popular algorithms based on the augmented Lagrangian technique are  {analytical target cascading} \cite{dormohammadi2012comparison},  {auxiliary problem principle} \cite{cohen1980auxiliary},  {dual decomposition method} \cite{2200000016}, and {proximal message passing} \cite{kraning2013dynamic}. One can also move away from the augmented Lagrangian, and use  {optimality condition decomposition} \cite{conejo2002decomposition},  {consensus methods}, or distributed algorithms developed from  {subgradient methods} \cite{4749425}, and  {dynamic programming} \cite{bertsekas1982distributed}.

\section{Online Convex Optimization} 

\subsection{Online optimization} The online optimization paradigm applies a robust optimization perspective for sequentially arriving data and costs that are  too complex to be efficiently modeled.  With its roots in classical ideas of sequential decisions in multi-armed bandit problems from the 1930s, online optimization has recently emerged as a prominent tool in machine learning, solving problems ranging from \emph{recommender systems} to \emph{spam filtering} \cite{Orabona2019AMI,OPT-013}. Online optimization portrays decisions for optimizing time-varying cost functions  as a feedback process, where one learns from experience as time evolves. Performance is considered  with respect to a static optimal decision taken in hindsight. Formally, the procedure of online optimization may be described as a game between a learner and an adversary played across a finite time horizon $t=1,\dots,T$.

\begin{mybox}{Online Optimization Paradigm \cite{Orabona2019AMI}}
For $t=1,\dots,T$, DO
\begin{itemize}
\item[(i)] The adversary selects a cost function $\ell_t(\cdot):\mathcal{X}\subseteq \mathbb{R}^d\to \mathbb{R}$ and keeps it to itself;
\item[(ii)] The learner makes a decision $\mathbf{x}_t\in \mathcal{X}$;
\item[(iii)] The learner suffers a loss $\ell_t(\mathbf{x}_t)$, and receives the cost function $\ell_t(\cdot)$ (full information), or just value of the loss  $\ell_t(\mathbf{x}_t)$ (bandit information). 
\end{itemize}
\tcblower
\end{mybox}

In sharp contrast to the view of classical optimization (i.e. a classical learner), where the loss function $\ell_t(\cdot)$ is revealed before the learner attempts to minimize it, online optimization acknowledges the difficulty in knowing  $\ell_t(\cdot)$ or even a model of it before decisions are made. The information that the learner receives about $\ell_t(\cdot)$ may be the whole function, a scenario referred to as {\em full information}; or the learner only experiences   losses at selected decisions, and in this case, we talk about {\em bandit information}. The loss functions $\ell_t(\cdot)$ are generally assumed to be arbitrary (but chosen from a given function class). Hence, it is impossible for the learner to infer $\ell_t(\cdot)$ before the decisions are made. As a result, it is sensible for the learner to identify $\mathbf{x}_1,\dots,\mathbf{x}_T \in \mathcal{X}$ so that {\em regret}, i.e.,
\begin{align*} 
\Reg(T):= \sum_{t=1}^T \ell_t (\mathbf{x}_t)- \min_{\mathbf{x}\in \mathcal{X}} \sum_{t=1}^T \ell_t (\mathbf{x})
\end{align*}
is minimized. From the definition, $\min_{\mathbf{x}\in \mathcal{X}} \sum_{t=1}^T \ell_t (\mathbf{x})
$ is the minimal accumulative loss of an oracle making a static decision to whom all $\ell_t(\cdot)$ are known   before $t=1$. Therefore, $\Reg(T)$ represents the difference between  the actual accumulative loss experienced by the learner compared to that of such an oracle, i.e., the  {regret}.

\subsection{Impact of feedback}

Let $\mathcal{X}\subseteq \mathbb{R}^d$ be a compact convex set containing the origin, for which $\mathbb{P}_{\mathcal{X}}$ is the projection onto $\mathcal{X}$. A simple yet effective  algorithm for the online learner is gradient descent implemented sequentially. The standard online gradient descent algorithm for solving the online optimization problem with full information is described below where $\alpha_t$ is the stepsize.

\begin{mybox}{Online Gradient Descent: Full Information Feedback \cite{OPT-013}}
 
For $t=1,\dots,T$, DO
\begin{itemize}
\item[(i)] The adversary selects a cost function $\ell_t(\cdot):\mathcal{X}\subseteq \mathbb{R}^d\to \mathbb{R}$ and keeps it to itself;
\item[(ii)] The learner makes a decision 
\begin{align}\label{eq:onlinegradient}
\mathbf{x}_t =\mathbb{P}_{\mathcal{X}}\Big(\mathbf{x}_{t-1} - \alpha_t \nabla \ell_{t-1}(\mathbf{x}_{t-1}) \Big);
\end{align}
\item[(iii)] The learner suffers a loss $\ell_t(\mathbf{x}_t)$, and receives  $\ell_t(\cdot)$. 
\end{itemize}
\tcblower
\end{mybox}

With bandit information, the learner only experiences losses and the loss function $\ell_t(\cdot)$ (and its gradient) is still unknown. Denote $\mathcal{K}_{\delta}=\left\{\mathbf{x}: \frac{1}{1-\delta} \, \mathbf{x} \in \mathcal{K}\right\}$.  Let $\mathbb{S}$ be the unit sphere in $\mathbb{R}^d$ under standard Euclidean norm. Then one can build unbiased gradient estimates from experienced losses to replace the true gradients in online gradient descent, leading to the following online bandit optimization algorithm. 
\begin{mybox}{Online Bandit Optimization: Bandit Information Feedback \cite{10.5555,OPT-013}}
Initialize $\mathbf{y}_0=0$. For $t=1,\dots,T$, the learner DO
\begin{itemize}
\item[(i)] Draw $\mathbf{u}_t \in \mathbb{S}$ uniformly at random.
\item[(ii)] Play $\mathbf{x}_t=\mathbf{y}_t +\delta \mathbf{u}_t$; receive loss $\ell_t(\mathbf{x}_t)$.
\item[(iii)] Build gradient estimate $\mathbf{g}_t= d \ell_t(\mathbf{x}_t) \mathbf{u}_t/\delta$; Update \begin{align}\label{eq:onlinebandit}
\mathbf{y}_{t+1} =\mathbb{P}_{\mathcal{X}_{\delta}}\big(\mathbf{x}_{t} - \eta \mathbf{g}_t \big).
\end{align}
\end{itemize}
\tcblower
\end{mybox}

It is immediately clear that in both (\ref{eq:onlinegradient}) and (\ref{eq:onlinebandit}), feedback is taking place. The promise of these online optimization algorithms lies in the fact that, when the stepsizes are selected as some suitable learning rates, the algorithms will produce   sub-linear regrets\footnote{For bandit feedback, the regret is technically $\mathbb{E} \Reg(T)$ where the expectation is taken over the randomness in the gradient estimate.} 
$\lim_{T\to \infty} \Reg(T)/T=0$,
for a suitably regular classes of convex cost functions. 
 This is a strong testimony to the performance of true {\em learning} during the sequential decisions.  The  regret averaged over time is  close to zero for sufficiently long time horizon: it is as if all the $\ell(\cdot)$ are known before the whole play starts and the learner decides to play a static optimal decision. With careful classification of the function classes for the loss functions, refined upper bounds on   $\Reg(T)$ can be established at $O (\log T)$, $O(\sqrt{T})$, $O (T^{3/4})$, etc \cite{OPT-013}. 
\subsection{Distributed online optimization}

In practice, the loss $\ell_t(\cdot)$ might represent a system-level loss,    scattered across a number of subsystems indexed in $\calV = \{1,\ldots,N \}$, such that $\ell_t(\cdot)=\sum_{i=1}^N\ell_{i,t}(\cdot)$. The overall system forms a network, where a directed graph $\mathcal{G}=(\mathcal{V},\mathcal{E})$ describes the communication structure of the network. Then the question arises on whether in the case that   subsystems may only talk to their neighbors over the graph $\mathcal{G}$, this will enable distributed online learning throughout the network. 
Following the distributed optimization framework, a distributed online optimization paradigm can be described as follows \cite{TAC} (see also \cite{pmlr-v70-zhang17g,7399359,9216151,8950383}). 

\begin{mybox}{Distributed Online Convex Optimization  \cite{TAC}}

{\bf Initialize} $\calX$ as a convex subset of $\reals^d$ defined by a family of inequalities: $\calX = \{ \bfx\in\reals^{d} \;|\; c_s(\bfx)\leq 0,\ s = 1,\ldots,p\}$.

\medskip

For $t=1,\dots,T$, agents in $\mathcal{V}$ DO
\begin{itemize}
\item Each agent $i\in \calV$ selects $\bfx_i(t)\in\mathcal{X}$, and a local adversary chooses $\ell_{i,t}(\cdot): \mathbb{R}^d \to \mathbb{R}$ as a convex cost function;
\item Each agent experiences a loss $\ell_{i,t}(\bfx_i(t))$;
\item The function $\ell_{i,t}$ is revealed to agent $i$; The decisions of the neighbors of the agent $i$ are also revealed to $i$, i.e., $\bfx_{j}(t)$ for $j\in \mathcal{N}_i:=\{(j,i)\in\mathcal{E}\}$.
\end{itemize}
\tcblower
\end{mybox}

The decision set $\mathcal{X}$ implies that in each time-step  each agent should be able  to perform a projection onto $\calX$, which can be   computationally expensive. Instead, one can only  require  that the constraints are satisfied in the long run, i.e., that $\sum_{t=1}^{T} \sum_{i=1}^{N} \sum_{s=1}^{p} c_s(\bfx_{i}(t))  \leq 0$. An effective  distributed online learning algorithm, then should aim to minimize the accumulated system-wide loss. The system-level regret is defined as the worst possible regret for all agents:
\begin{equation}
\begin{aligned}
\Sreg(T)
&:= \max_{i\in \mathrm{V}} \left[ \sum_{t=1}^{T} \sum_{j=1}^{N} \ell_{j,t}(\bfx_{i}(t)) - \sum_{t=1}^{T} \sum_{j=1}^{N} \ell_{j,t}(\bfx^{\star})\right ]
\label{regret}
\end{aligned}
\hspace{10000pt minus 1fil} 
\end{equation}
where $\bfx^\star = \arg\min_{\bfx\in\calX} \sum_{t=1}^{T} \sum_{j=1}^{N} \ell_{j,t}(\bfx)$ is the system-level decision by a static optimal oracle. The performance of the algorithm is further characterized by the so-called cumulative {absolute} constraint violation defined by
\begin{equation}
\begin{aligned}
\Cacv(T) := \sum_{t=1}^{T}\sum_{i=1}^{N} \sum_{s=1}^{p} [c_s(\bfx_{i}(t))]_{+}
\label{CACV}
\end{aligned}
\end{equation}
where  $\left[a\right]_+ =\max\{ 0,a\}$.

\subsection{Distributed online primal-dual gradient algorithm}
In classical distributed optimization, it is popular to use a consensus algorithm as an information aggregation subroutine. Specifically, we may associate  a doubly stochastic matrix $\mathbf{A}$ with the graph $\mathcal{G}=(\mathcal{V},\mathcal{E})$ such that $\mathbf{A}_{ij}>0$ if and only if $(j,i)\in \mathcal{E}$. In general, for strongly connected graph $\mathcal{G}$, one can always find such an $\mathbf{A}$. 

\vspace{2mm}
\noindent{\em Example 3.}
We illustrate the performance of the proposed algorithms using a simple experiment. Specifically, we consider a distributed online linear regression problem over a network, where $\ell_{i,t}(\bfx)=\left( \mathbf{a}_i(t)^\top \bfx - b_i(t) \right)^2/2 $. The constraints are described as 
\begin{align}
 c_m(\bfx) &= L - \bfx_m \leq 0, \quad m=1,\ldots,d,  \\
  c_{d+m}(\bfx) &= \bfx_m - U \leq 0, \quad m=1,\ldots,d. 
\end{align}
Every entry of $\mathbf{a}_i(t)$ and $b_i(t)$ is generated uniformly at random within the interval $[-1,1]$ and $[0,1]$, respectively, independently for each time $t=1,\dots,T$. Throughout the experiments, we implement the distributed online primal-dual gradient algorithms proposed in \cite{TAC} with full information or bandit information feedback.

\noindent{\em System Setup}. The graph $\mathrm{G}$ is randomly generated and selected  as depicted in Fig. 1. The weighting matrix associated with the network in Fig. 1 is generated according to the maximum-degree weights:
\begin{eqnarray}
\bfA_{ij}
&=&
\left\{
\begin{array}{ll}
\frac{1}{1 + d_{\max}},           &\qquad (j,i)\in\calE  \\
1 - \frac{d_i}{1 + d_{\max}} ,    &\qquad i=j \\
0,                                &\qquad (j,i)\notin\calE
\end{array}
\right.
\label{weight-max}
\end{eqnarray}
where $d_{\max} = \max_{i\in\mathrm{V}} \{d_i\}$ is the maximum degree of $\mathcal{G}$ with $d_i$ denoting the degree of node $i$. We set the parameters as follows: $N = 20$, $d = 2$, $L = -1/2$, $U = 1/2$, and $R_\calX = U  \sqrt{d}$. The performance of the algorithm  is averaged over 10 runs.

\begin{figure}[H]
  \centering
  \includegraphics[width=.9\linewidth]{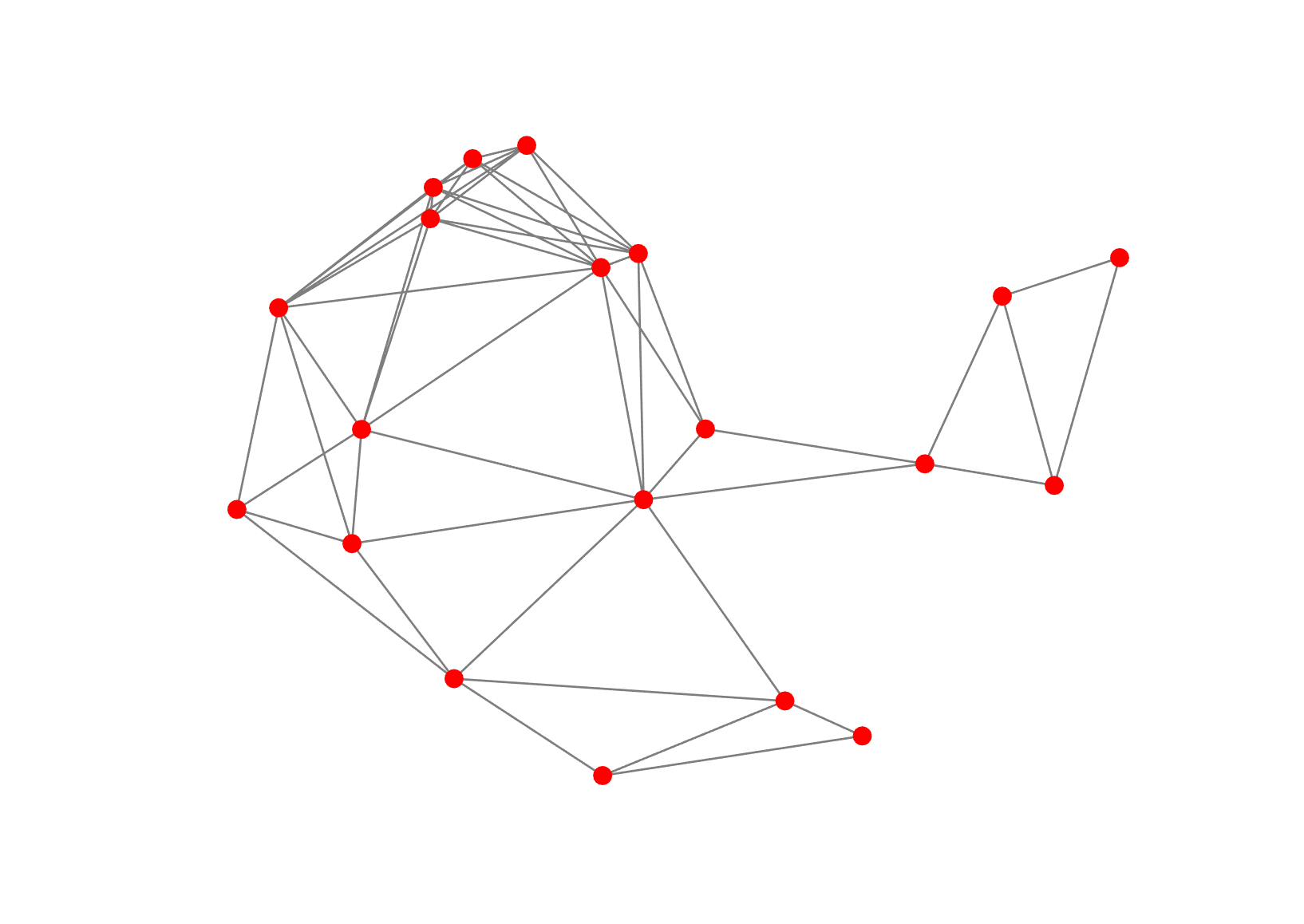}
  \caption{A randomly generated network of 20 nodes.}
\end{figure}

\begin{figure}[H]
  \centering
  \includegraphics[width=1.1\linewidth]{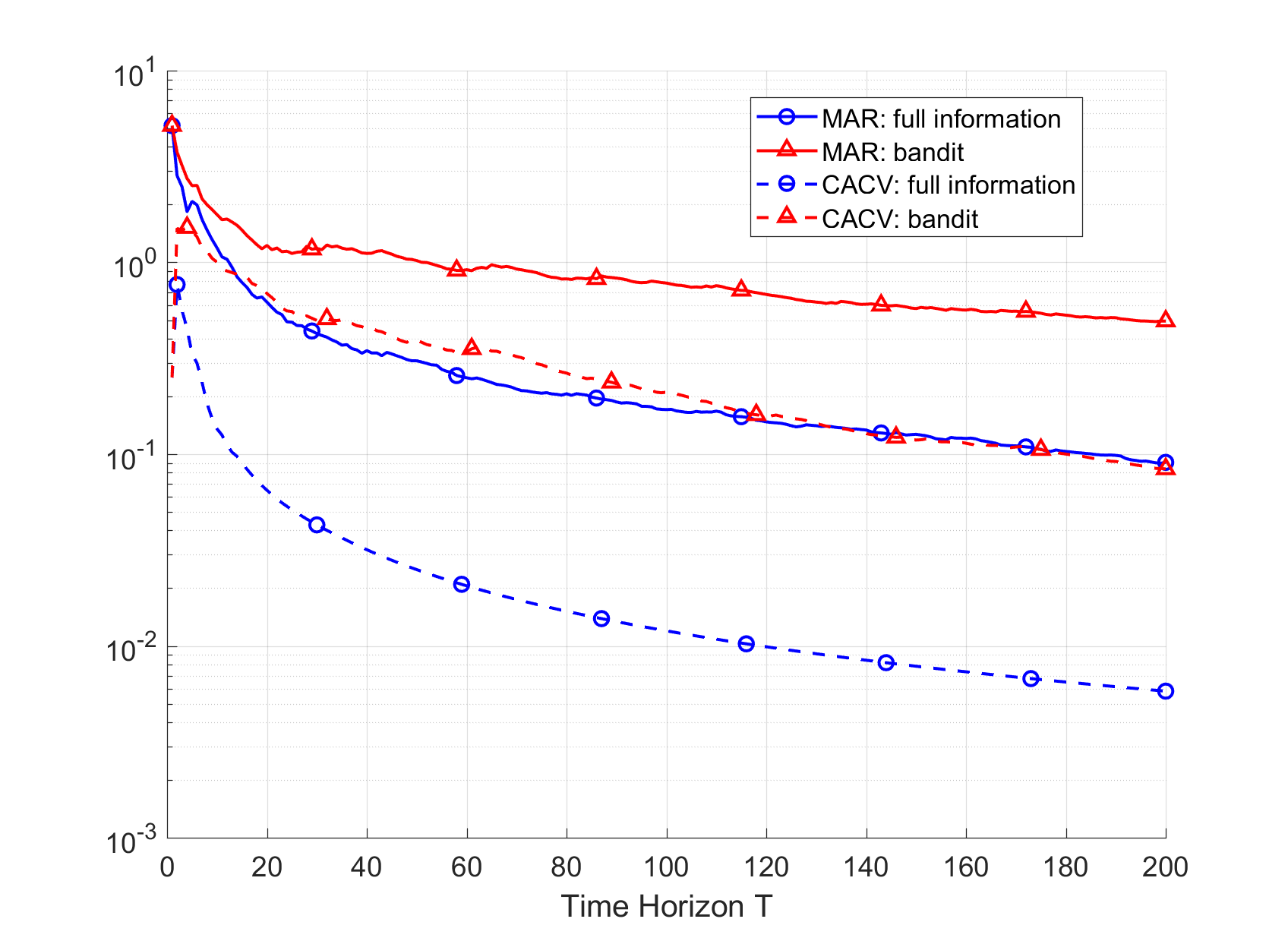}
  \caption{ASR and ACV vs. time for the distributed online optimization algorithms with full information and bandit information  in \cite{TAC}.}\label{figperformance}
\end{figure}

\noindent{\em Performance.} We run the algorithms and plot the Average System Regret (ASR, for short) defined as {\sf SReg}$(T)/T$ and the Average Constraint Violations (ACV, for short)  defined as $\Cacv(T)/T$, as a function of the time horizon $T$ in Fig. \ref{figperformance}.  Clearly both the ASR and ACV converge to zero. 

\begin{figure}[H]
  \centering
  \includegraphics[width=1.1\linewidth]{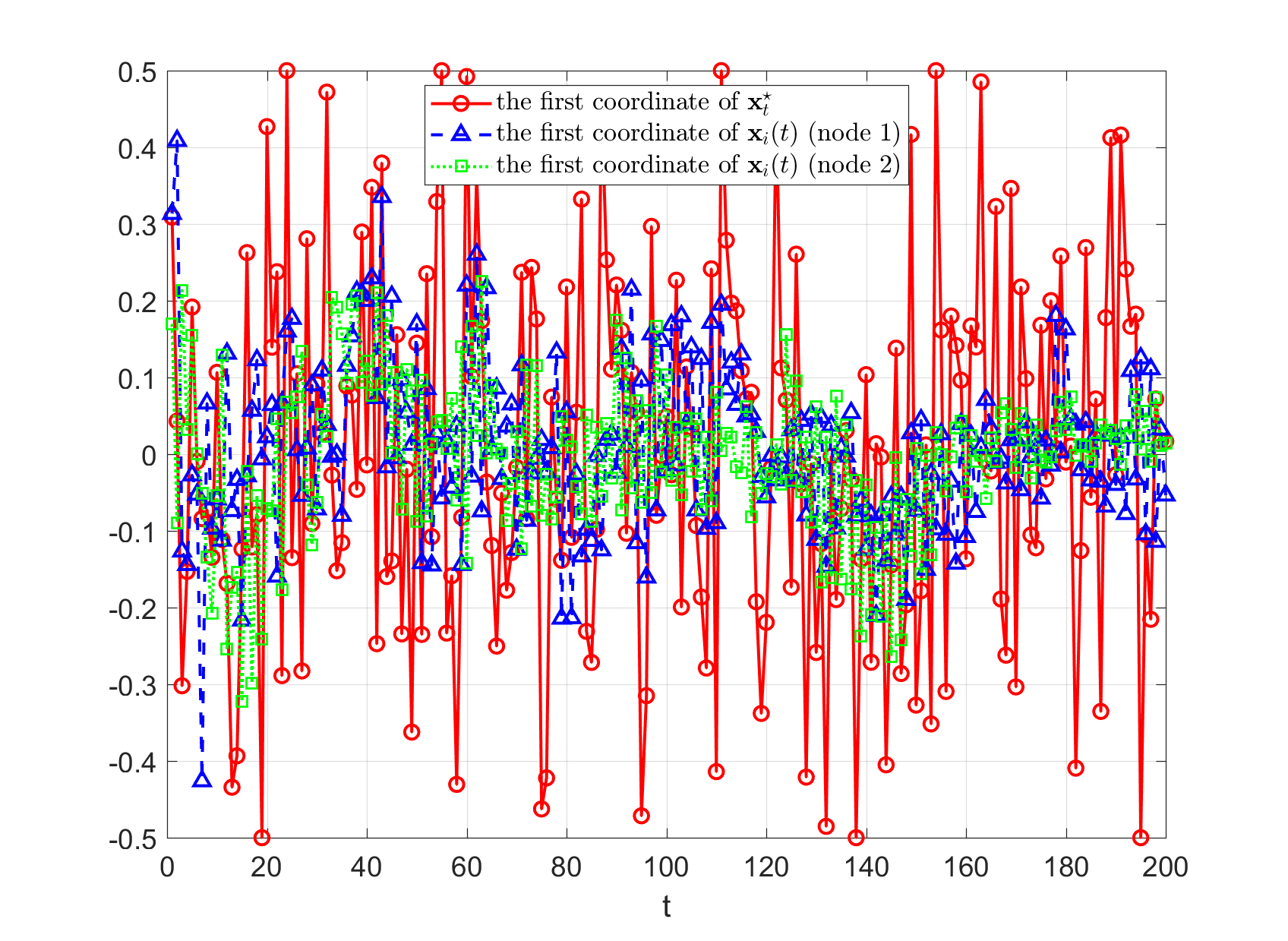}
  \caption{System-level optimal decisions from repeated offline optimization vs. distributed agent decisions from sequential online optimization.}\label{figstationary}
\end{figure}

\noindent{\em Decision Stationarity.} At time $t$, let the system-level objective function 
$\sum_{i=1}^{N} \frac{1}{2} \left( \mathbf{a}_i(t)^\top \bfx - b_i(t) \right)^2$ yield an optimal decision $\bfx_t^\ast$. In Fig. \ref{figstationary}, we plot the first entry of the repeated offline optimizer $\bfx_t^\ast$, and the first entry of the distributed sequential online optimizer $\mathbf{x}_1(t)$ and  $\mathbf{x}_2(t)$ for agent 1 and agent 2, respectively. It can be seen that the online agent decisions demonstrate significantly reduced fluctuations compared to the repeated system-level offline decisions. 

\section{Perspectives on Online Optimization for Power Grid}

\subsection{Offline vs. online optimization}
For problems such as the economic dispatch and optimal power flow in Example 1 and Example 2 over a time horizon $t=1,\dots,T$, there may be two paradigms. 

\vspace{2mm}

\noindent{\em Distributed (Repeated Offline) Optimization [DRO-O]}. For each time $t$,  independently treat the corresponding problem (\ref{ex1}) and (\ref{ex2}); apply distributed optimization algorithms until suitable convergence is guaranteed for time $t$; implement the optimal decision $\mathbf{x}_i^\ast(t)$ for $i=1,\dots,N$ repeatedly at the respective time $t$. 

\vspace{2mm}

\noindent{\em Distributed Online Optimization [DO-O]}. Employ the Distributed Online Convex Optimization paradigm outlined in Section 2.3; apply distributed online optimization algorithms   throughout the horizon $t=1,\dots,T$;  implement the   decision $\mathbf{x}_i (t)$ for $i=1,\dots,N$ sequentially for $t=1,\dots,T$. 

\vspace{2mm}

The potential in developing online optimization frameworks for problems in power grids has drawn attentions in the literature.   Online convex optimization has been adopted  in \cite{zhou2017incentive} for the control of distributed energy sources in the context of social welfare maximization. Under a similar social welfare maximization paradigm, \cite{zhou2019online} considers control of distributed energy sources with both continuous and discrete constraints. Moreover, \cite{li2019distributed} provides a unified framework for economic dispatch and unit commitment and proposes a centralized and distributed online convex optimization method for exploring such a framework. 

Next, we would like to offer a few perspectives towards the strengths, challenges, and possible future direction for online optimization in power grids.

\subsection{Perspectives between [DRO-O] and   [DO-O] }
First, it is worth mentioning that the key difference between [DRO-O] and [DO-O] goes far beyond the respective classes of algorithms. Underpinning the two frameworks are fundamentally different views about the system: 
\begin{itemize}
\item In [DRO-O], the time-varying cost functions are {\em known before  decisions}; 
\item In [DO-O], the time-varying cost functions are {\em experienced after decisions}.  
\end{itemize}
As a result, conceptually the [DRO-O] algorithms are optimizers, while the [DO-O] algorithms are learners. Therefore, [DO-O] suits systems that are uncertain or unpredictable. 

Next, the strength of [DO-O] lies in guaranteed sub-linear regret against adversaries. In practice, the adversaries represent the worst-case scenarios. Remarkably, the aforementioned regret bounds of orders $O (\log T), O(\sqrt{T}), O (T^{3/4})$ may be valid even for feedback adversaries, where the cost function $\ell_t(\cdot)$ depends on the past experiences. Moreover, the online decisions in [DO-O]  tend to converge to a static optimal decision with respect to the cumulative cost over the entire horizon, while decisions in [DRO-O] tend not to converge as they are tracking real-time optimal decisions for time-varying cost functions. This is shown in Example 3 where online decisions indeed are more stationary compared to repeated offline optimal decisions. 

In the context of power grids, [DO-O] might be more suited in problems related to wind or solar energy grid-integration, and energy storage applications including residential batteries and EVs. Such applications involve significant uncertainty regarding the weather, network impedance and topology, real-time price volatility, and user preferences including when, where and for how long an EV will require charging. Importantly, for problems with known grid and user information, [DRO-O] is a more sensible choice as the performance of [DO-O] is much more conservative.

\subsection{Future directions}
Towards establishing practical online optimization frameworks  for problems in power grid, there are a few possible  directions. First, the notion of regret needs to be taken into account for online optimization of power grid problems. Existing regret bounds for online optimization are for classes of convex, smooth, or strongly convex functions, etc. Cost functions in power grid problems and constraints are certainly more structured (despite being unknown before decisions), and thus refined regret bounds might exist. Second, the uncertain nature of online optimization needs to be carefully matched to practice. The characterization of cost functions  should also be evaluated in the power system context. Third, hybrid decision frameworks that combine the strength of [DRO-O] and [DO-O], where the information and uncertainty of the cost functions can be jointly treated, would be of significant value for power grid applications.



\begin{acks}
We wish to thank Professor Alexandre Proutiere for insightful feedback and suggestions. This work was supported by the Australian Research Council under Grants DP180101805 and DP190103615. 
\end{acks}


\appendix

\end{document}